\documentclass[11pt,a4paper,leqno,twoside, headinclude]{amsart}

\title{Computational complexity of topological invariants}
\def\titl{Computational complexity of topological invariants}
\def\auth{Manuel Amann}
\date{December 4th, 2011}

\subjclass[2010]{68Q17 (Primary), 55P62 (Secondary)}
\keywords{\noindent computational complexity, cup-length, LS-category, pure elliptic space}
\thanks{The author was supported by a Research Grant of the German Research Foundation.}
\author{\auth}

\usepackage[english]{babel}

\usepackage[colorlinks,pdftex, plainpages=false]{hyperref}
\hypersetup{pdftitle=\titl, pdfauthor=\auth, pdftoolbar=false,
plainpages=false, hyperindex=true, pdfdisplaydoctitle=true}

\hypersetup{colorlinks,%
citecolor=black,%
filecolor=black,%
linkcolor=black,%
urlcolor=black,%
pdftex}

\usepackage{color}

\usepackage{amsmath, amssymb, amscd, amsthm}
\usepackage{stmaryrd}
\usepackage{latexsym}
\usepackage[all]{xy}
\usepackage{pb-diagram}
\usepackage{rotating}
\usepackage{multicol}
\usepackage{lscape}
\usepackage{wasysym}
\usepackage{longtable}
\usepackage{enumerate}
\xyoption{all}

\newtheorem{theo}{Theorem}[section]
\newtheorem{main}{Theorem}
\newtheorem{maincor}[main]{Corollary}
\newtheorem*{main*}{Theorem}

\newtheorem*{mainprop*}{Proposition}

\newtheorem{mainconj}{Conjecture}

\newtheorem{prop}[theo]{Proposition}
\newtheorem{defi2}[theo]{Definition}
\newtheorem*{defi2*}{Definition}

\newenvironment{defi*}{\begin{defi2*}\normalfont}{\end{defi2*}}
\newenvironment{defin}[1]{\begin{defi2}[#1]\normalfont}{\end{defi2}}
\newenvironment{defin*}[1]{\begin{defi2*}[#1]\normalfont}{\end{defi2*}}
\newtheorem{rem2}[theo]{Remark}
\newenvironment{rem}{\begin{rem2}\normalfont}{\hfill$\boxbox$\end{rem2}}

\newtheorem{lemma}[theo]{Lemma}

\newtheorem*{cor*}{Corollary}

\newtheorem*{conj*}{Conjecture}
\newtheorem*{theo*}{Theorem}
\newtheorem*{ques*}{Question}
\newtheorem*{mi2}{Main Idea}

\newtheorem{ex2}[theo]{Example}

\newtheorem{exer2}[theo]{Exercise}

\newtheorem{alg2}[theo]{Algorithm}

\newcommand{\cc}{{\mathbb{C}}}                                     % complex numbers
\newcommand{\nn}{{\mathbb{N}}}                                     % natural numbers
\newcommand{\qq}{{\mathbb{Q}}}                                     % rational numbers
\newcommand{\zz}{{\mathbb{Z}}}                                     % integers
\newcommand{\dif} {{\operatorname{d}}}                             % differential operator d
\newcommand{\In} {{\,\subseteq\,}}                                 % subset
\newcommand{\im} {{\operatorname{im\,}}}                           % image
\newcommand{\id}{{\operatorname{id}}}                              % identity
\newcommand{\APL}{{\operatorname{A_{PL}}}}                         % polynomial differential forms
\newcommand{\cat} {{\operatorname{cat}}}                           % Lusternik--Schnirelman category
\newcommand{\cupl}{{\operatorname{c}}}                             % Cup-length
\newcommand{\e}{{\operatorname{e}}}                                % Toomer invariant
\def\co{\colon\thinspace}                                          % colon in maps

\newcommand{\comment}[1]{}                                         % insert a large comment
\newcommand{\hto}[1]{\overset{#1}{\hookrightarrow}}                % abbreviation for labelled \hookrightarrow
\newcommand{\str}{\noindent\textbf{Structure of the article. }}    % Structure of the article
\newcommand{\odd}{\textrm{odd}}                                    % odd
\newcommand{\even}{\textrm{even}}                                  % even

\newenvironment{prf}{\begin{proof}[\textsc{Proof}]} {\end{proof}}     % alternative environment proof
\begin{document}

\maketitle \thispagestyle{empty}

%%%%%%%%%%%%%%%%%%%%%%%%%%%%%%%%%% Abstract%%%%%%%%%%%%%%%%%%%%%%%%%%%%%%%%%%%%%%%%

\begin{abstract}
 We answer the following question posed by Lechuga: Given a simply-connected space $X$ with both $H_*(X,\qq)$ and $\pi_*(X)\otimes \qq$ being finite-dimensional, what is the computational complexity of an algorithm computing the cup-length and the rational Lusternik--Schnirelmann category of $X$?

Basically, by a reduction from the decision problem whether a given graph is $k$-colourable (for $k\geq 3$) we show that (even stricter versions of the) problems above are $\mathbf{NP}$-hard.
\end{abstract}

%%%%%%%%%%%%%%%%%%%%%%%%%%%%%%%%%% Introduction %%%%%%%%%%%%%%%%%%%%%%%%%%%%%%%%%%%

\section*{Introduction}

The theory of computational complexity has developed a powerful machinery of describing how ``difficult'', i.e.~how time-consuming, it is to answer certain posed questions algorithmically. Most classically, this asks for the following categorification of problems: The complexity class $\mathcal{P}$ describes all the problems for which there is a polynomial-time solving algorithm; the class $\mathbf{NP}$ is formed by those problems which may at least be verified in polynomial time. Clearly $\mathbf{P}\In \mathbf{NP}$, however, it is the common belief that several problems in $\mathbf{NP}$ are much harder to solve than the problems in $\mathbf{P}$. Known algorithms typically run at exponential costs.

A whole variety of problems stemming from completely different areas of mathematics and computer science have been found to be harder than all the problems in $\mathbf{NP}$, i.e.~to be $\mathbf{NP}$-hard. Just to name a few most prominent ones we mention the \emph{knapsack problem} and the \emph{subset sum problem}, the \emph{Hamilton circuit problem} and the \emph{travelling salesman problem}, the \emph{satisfiability problem} and the \emph{graph colouring problem}.

Also in the field of algebraic topology it is easy to imagine several problems for which it seems difficult to find efficient solving algorithms. In particular, Rational Homotopy Theory has the appeal of providing ``computable problems'', which certainly ask for algorithmic treatment. Indeed, Rational Homotopy Theory permits a categorical translation from topology/homotopy theory to algebra at the expense of losing torsion information. Yet, it turns out that the algebraic side allows for concrete computations.

Using this approach several topological problems were shown to be $\mathbf{NP}$-hard. In \cite{Ani89} it is shown that computing the rational homotopy groups $\pi_*(X)\otimes \qq$ of a simply-connected CW-complex $X$ is $\mathbf{NP}$-hard. So is the problem of whether a simply-connected space $X$ with $\dim \pi_*(X)\otimes \qq<\infty$ also has finite-dimensional rational cohomology (cf.~\cite{GL03}). In the same article it was shown that for formal spaces, i.e.~for spaces for which the rational homotopy type can be formally derived from the rational cohomology algebra, the computation of Betti numbers, of cup-length and of the rational Lusternik--Schnirelmann category are $\mathbf{NP}$-hard problems. In \cite{GL03} it is shown that the computation of Betti numbers of a simply-connected space with both finite-dimensional rational homotopy and finite-dimensional rational homotopy, a \emph{(rationally) elliptic space}, is $\mathbf{NP}$-hard.

However, already in the article \cite{Lec02} and then explicitly in \cite{GL03} the following question is posed:
\begin{ques*}[Lechuga]
Given an elliptic space, what is the computational complexity of computing its rational cup-length or its rational Lusternik--Schnirelmann category?
\end{ques*}

Note that the methods from the results presented above do not answer this question, as they do not apply to the case of elliptic spaces. (For a definition of the topological invariants we refer the reader to section \ref{sec04}.)

In this article we shall answer Lechuga's question by revealing these problems as $\mathbf{NP}$-hard---we do so already for the same question posed on the subclass of pure elliptic spaces. For this we specify the following problems

\vspace{3mm}

\begin{itemize}
\item[$\mathcal{P}$:]
Let $X$ be a simply-connected topological space with finite-dimensio\-nal rational homology $H_*(X,\qq)$ and with finite-dimensional rational homotopy $\pi_*(X)\otimes \qq$. What is its cup-length?
\end{itemize}

\vspace{3mm}

\begin{itemize}
\item[$\mathcal{Q}$:]
Let $X$ be a simply-connected topological with finite-dimensional rational homology $H_*(X,\qq)$ and with finite-dimensional rational homotopy $\pi_*(X)\otimes \qq$. What is its rational Lusternik--Schnirelmann category?
\end{itemize}

\vspace{3mm}

We are interested in the computational complexity of the problems $\mathcal{P}$ and $\mathcal{Q}$.

\vspace{5mm}

The codification of a simply-connected space $X$ will be given as the data contained in its \emph{minimal Sullivan model} $(\Lambda V_X,\dif)$, i.e.~$X$ will be represented by the degrees of the homogeneous generators $x_1,\dots, x_l$ of $V_X$ and the coefficients of the polynomials in the $x_i$ which represent the differential.

Recall that a minimal Sullivan model of a simply-connected space is a free (graded) commutative graded algebra $\Lambda V$ over the $\zz$-graded rational vector space $V=V^{\geq 2}$ together with a differential $\dif$ defined by $\dif\co V^*\to (\Lambda V)^{*+1}$ and extended to $\Lambda V$ as a derivation. The differential satisfies that its image lies in the subalgebra of elements of wordlength at least two in $V$, i.e.~$\im \dif \In \Lambda^{\geq 2} V$. One then requires the existence of a quasi-isomorphism $(\Lambda V,\dif)\to \APL(X)$, i.e.~a morphism of differential graded algebras to the polynomial differential forms $\APL(X)$ on $X$ inducing an isomorphism on homology. Thus $(\Lambda V,\dif)$ encodes the rational homotopy type of $X$. (See \cite{FHT01}.3, \cite{FHT01}.10 and \cite{FHT01}.12 for the missing definitions.)
In particular, the homology algebra $H(\Lambda V,\dif)$ of the minimal model is the cohomology algebra of $X$.

\vspace{5mm}

Thus problem $\mathcal{P}$ translates to problem

\vspace{3mm}

\begin{itemize}
\item[$\mathcal{P'}$:]
Let $(\Lambda V,\dif)$ be a simply-connected elliptic Sullivan algebra. What is its cup-length?
\end{itemize}

\vspace{3mm}

Since the rational Lusternik--Schnirelmann category of a simply-connected space with rational homology of finite type equals the category of its minimal Sullivan model (cf.~\cite{FHT01}.29.4, p.~386), problem $\mathcal{Q}$ becomes

\vspace{3mm}

\begin{itemize}
\item[$\mathcal{Q}'$:]
Let $(\Lambda V,\dif)$ be a simply-connected elliptic Sullivan algebra. What is its Lusternik--Schnirelmann category?
\end{itemize}

\vspace{3mm}

By the same theorem the rational Toomer invariant $\e_0(X)$ of a simply-connected space with rational homology of finite type equals the Toomer invariant $e(\Lambda V,\dif)$ of its minimal model. The rational cohomology algebra of a simply-connected elliptic space satisfies Poincar\'e duality. On simply-connected spaces with cohomology satisfying Poincar'e duality the rational Toomer invariant equals the rational Lusternik--Schnirelmann category (cf.~\cite{FHT01}.38, p.~511). Thus problems $\mathcal{Q}$ and $\mathcal{Q}'$ have the obvious analogue

\vspace{3mm}

\begin{itemize}
\item[$\tilde{\mathcal{Q}}$:]
Let $(\Lambda V,\dif)$ be a simply-connected elliptic Sullivan algebra; (respectively let $X$ be a simply-connected elliptic space). What is its (rational) Toomer invariant?
\end{itemize}

\vspace{3mm}

A Sullivan algebra $(\Lambda V,\dif)$ is \emph{pure} if $V=P\oplus Q$ with $Q=V^\even$ and $P=V^\odd$ and the differential $\dif$ satisfies
\begin{align*}
\dif|_Q=0 \qquad \textrm{and } \qquad \dif(P)\in \Lambda Q
\end{align*}
Classical examples of spaces admitting pure models are biquotients; respectively, in particular, their subclass of homogeneous spaces.

\vspace{5mm}

We shall determine the computational complexities of the following stricter problems, i.e.~we shall show that they are $\mathbf{NP}$-hard. In particular, this will answer the original question by Lechuga.

\vspace{3mm}

\begin{itemize}
\item[$\mathcal{P}''$:]
Let $(\Lambda V,\dif)$ be a simply-connected pure elliptic Sullivan algebra. What is its cup-length?
\end{itemize}

\vspace{3mm}

\begin{itemize}
\item[$\mathcal{Q}''$:]
Let $(\Lambda V,\dif)$ be a simply-connected pure elliptic Sullivan algebra. What is its rational Lusternik--Schnirelmann category?
\end{itemize}

\vspace{3mm}

This leads us to our main theorems.

\begin{main}\label{theoA}
The problem $\mathcal{P}''$ is $\mathbf{NP}$-hard.
\end{main}

\begin{main}\label{theoB}
The problem $\mathcal{Q}''$ is $\mathbf{NP}$-hard.
\end{main}

Since pure elliptic spaces form a subclass of elliptic spaces, we obtain
\begin{maincor}\label{corC}
The problems $\mathcal{P}$ and $\mathcal{P}'$ are $\mathbf{NP}$-hard. So are the problems $\mathcal{Q}$, $\mathcal{Q}'$ and $\tilde{\mathcal{Q}}$.
\end{maincor}

\vspace{3mm}

\str In section \ref{sec01} we briefly review some basic concepts from the theory of computational complexity.  We recall the definitions of the topological invariants in section \ref{sec04} before we prove theorem \ref{theoA} in section \ref{sec02}. Section \ref{sec03} is devoted to the proof of theorem \ref{theoB}.

%%%%%%%%%%%%%%%%%%%%%%%%%%%%%%%%%% Section 1 %%%%%%%%%%%%%%%%%%%%%%%%%%%%%%%%%%%%%%

\section{Basic notions from complexity theory}\label{sec01}

Let us recall some definitions from complexity theory.
\begin{defin}{problem, solution, complexity}\label{Gdef04}
A \emph{problem} is a set $\mathcal{X}$ of ordered pairs $(I,A)$ of bitcoded strings with the property that for each instance $I$ there exists an answer $A$ such that $(I,A)\in P$.

A \emph{decision problem} is a function $\mathcal{X}$ with values in $\{0,1\}$.

A \emph{solution} of a problem is an algorithm which computes for each input $I$ an output $A$ such that $(I,A)\in \mathcal{X}$ in finitely many steps.

The \emph{complexity} of a problem is the infimum of the (asymptotic) run times of all solution algorithms.
\end{defin}

\begin{defin}{$\mathbf{P}$, $\mathbf{NP}$}\label{Gdef05}
Suppose given an instance $I$ and a suggested proof $A_I$ for the fact that $(I,1)\in \mathcal{X}$ for the decision problem $\mathcal{X}$. A \emph{polynomial verifier} of $\mathcal{X}$ is an algorithm which checks in polynomial time whether $A_I$ really proves that $I$ is true.

The class of all decision problems for which there exists a polynomial time solution algorithm is called the class $\mathbf{P}$. The class of all decision problems for which there exists a polynomial verifier form the class $\mathbf{NP}$.

A decision problem $\mathcal{Y}$ is \emph{$\mathbf{NP}$-complete}, if all $\mathcal{X}\in \mathbf{NP}$ can be reduced to $\mathcal{Y}$ in polynomial time.

An arbitrary problem $\mathcal{X}$ which is harder than all the problems $\mathcal{Y}\in \mathbf{NP}$ is called \emph{$\mathbf{NP}$-hard}, i.e.~for each problem $\mathcal{Y}\in \mathbf{NP}$ an algorithm solving $\mathcal{X}$ can be translated in polynomial time to an algorithm solving $\mathcal{Y}$.
\end{defin}

In order to show that a problem is $\mathbf{NP}$-hard one tends to use a reduction principle: If one can reduce an $\mathbf{NP}$-complete problem $\mathcal{A}_1$ to a problem $\mathcal{A}_2$, then the latter has to be $\mathbf{NP}$-hard. Indeed, since $\mathcal{A}_1$ is $\mathbf{NP}$-complete, it is maximally hard in $\mathbf{NP}$, i.e.~every problem $\mathcal{A}\in \mathbf{NP}$ can be reduced to $\mathcal{A}_1$. Consequently, every problem $\mathcal{A}\in \mathbf{NP}$ can be reduced to $\mathcal{A}_2$ in polynomial time. Thus also $\mathcal{A}_2$ is harder than all the problems in $\mathbf{NP}$.

%%%%%%%%%%%%%%%%%%%%%%%%%%%%%%%%%% Section 4 %%%%%%%%%%%%%%%%%%%%%%%%%%%%%%%%%%%%%%

\section{The topological invariants} \label{sec04}

In this section we intend to briefly recall the definitions of the topological invariants which appear in abundance in this article. They will partly be defined using Rational Homotopy Theory. We recommend the textbook \cite{FHT01} for an introduction to this field. We shall follow the notation and definitions provided there.

Let us start with the simplest invariant which is provided by
\begin{defin}{(rational) cup-length}\label{def01}
The \emph{(rational) cup-length} $\cupl_0(X)$ of a path-connected topological space $X$ is the greatest number $n\in \zz \cup \{\infty\}$ such that there are cohomology classes $[x_1], \dots, [x_n] \in H^{> 0}(X,\qq)$ satisfying $[x_1]\cup \ldots \cup [x_n]\neq 0$, i.e.~
\begin{align*}
\cupl_0(X)=\sup \{n\in \nn_0 \mid \exists [x_1], \ldots, [x_n] \in H^{>0}(X,\qq):~ [x_1]\cup \ldots \cup [x_n]\neq 0\}
\end{align*}
\end{defin}

\vspace{5mm}

Let us now introduce the notion of Lusternik--Schnirelmann category and its rational analogue.

 A subset $U\In X$ of a topological space $X$ is called \emph{contractible in $X$} if its inclusion $i\co U\hto{} X$ is homotopic to a constant map.
\begin{defin}{(rational) Lusternik--Schnirelmann category of a space}\label{def02}
The \emph{Lusternik--Schnirelmann category} $\cat X$ of a topological space $X$ is the least number $m\in \zz\cup \{\infty\}$ such that $X$ is the union of $m+1$ open subsets $U_i$, each contractible in $X$.

The \emph{rational Lusternik--Schnirelmann category} of $X$ is the least number $m\in \zz\cup \{\infty\}$ such that $X\simeq_\qq Y$ and $\cat Y=m$.
\end{defin}
We shall mainly draw on the definition of category in the setting of Sullivan algebras. In order to provide a definition in this case we suppose that $(\Lambda V,\dif)$ is a Sullivan algebra and that $m\geq 1$.

Taking the quotient of $(\Lambda V,\dif)$ by all the elements $\Lambda^{> m} V$ of wordlength larger than $m$ induces the structure of a commutative cochain algebra for $(\Lambda V/\Lambda^{>m} V,\dif)$. This is due to the fact that the differential $\dif$ is a derivation.
The surjection
\begin{align*}
f_m\co (\Lambda V,\dif)\to (\Lambda V/\Lambda V^{>m},\dif)
\end{align*}
extends to a model $\varphi_m$ of $f_m$.
\begin{align}\label{eqn01}
\xymatrix{
(\Lambda V,\dif) \ar[dr]_{f_m} \ar@{^{(}->}[r]^>>>>>{i_m}& {(\Lambda V\otimes \Lambda Z(m),\dif)}
\ar[d]^{\varphi_m}_\simeq\\
&(\Lambda V/\Lambda^{>m} V,\dif)
}
\end{align}

With this notation we make
\begin{defin}{Lusternik--Schnirelmann category of an algebra}\label{def03}
The \linebreak[4]\emph{Lusternik--Schnirelmann category} $\cat(\Lambda V,\dif)$ of a Sullivan algebra $(\Lambda V,\dif)$ is the least number $m\in \zz\cup \{\infty\}$ such that there is a cochain algebra morphism
\begin{align*}
p_m\co (\Lambda V\otimes Z(m),\dif)\to (\Lambda V,\dif)
\end{align*}
such that $p_m\circ i_m=\id_{\Lambda V}$.
\end{defin}
The rational category of a simply-connected topological space with rational homology of finite type equals the category of a respective Sullivan model $(\Lambda V,\dif)$---cf.~proposition \cite{FHT01}.29.4, p.~386.

Clearly, one always has $\cat_0 X\leq \cat X$. For simply-connected CW-complexes one obtains that $\cat_0 X=\cat X_\qq$ (cf.~proposition \cite{FHT01}.28.(i), p.~371).

\vspace{5mm}

Let us eventually briefly comment on Toomer's invariant.

\begin{defin}{(rational) Toomer invariant}
The \emph{Toomer invariant} \linebreak[4] $\e(X,\qq)$ of a topological space is the least number $m$ for which there is a continuous map $f\co Z\to X$ from an $n$-cone $Z$ (cf.~the definition on \cite{FHT01}, p.~359) such that $H^*(f,\qq)$ is injective.

The \emph{rational Toomer invariant} $\e_0(X)$ is the least number $m$ such that $X\simeq_\qq Y$ and $\e(Y,\qq)=m$.

In the notation of diagram \ref{eqn01} we define the \emph{Toomer invariant of a Sullivan algebra} $\e(\Lambda V,\dif)$ to be the least value $m\in \zz \cup \{\infty\}$ such that $H(f_m,\qq)$ is injective.
\end{defin}
Again, the rational Toomer invariant of a simply-connected topological space with rational homology of finite type equals the Toomer invariant of a respective Sullivan model $(\Lambda V,\dif)$---cf.~proposition \cite{FHT01}.29.4, p.~386.

\vspace{5mm}

We remark that in the light of the cited results we may use results formulated for simply-connected topological spaces with rational homology of finite type and translate them to respective Sullivan models.

%%%%%%%%%%%%%%%%%%%%%%%%%%%%%%%%%% Section 2 %%%%%%%%%%%%%%%%%%%%%%%%%%%%%%%%%%%%%%

\section{Proof of theorem A}\label{sec02}

Let $G=(V,E)$ be a (non-directed) finite connected simple graph with vertices $V=\{v_1,\dots, v_n\}$ and edges $E=\{(e_i,e_j) \mid (i,j)\in J\}$ for some index set $J$. Following \cite{LM00}, p.~91, we associate to $G$ and a given integer $k\geq 2$ a finitely generated pure Sullivan algebra by
\begin{align*}
V^\even_{G,k}&=\langle x_1, \dots,x_n\rangle
\intertext{with $\deg x_i=2$  and $\dif x_i=0$ for all $1\leq i\leq n$, by}
V^\odd_{G,k}&=\langle y_{i,j} \rangle_{(i,j)\in J}
\end{align*}
with $\deg y_i=2k-3$ and by
\begin{align*}
\dif y_{i,j}=\sum_{l=1}^k x_r^{k-1}x_s^{l-1}
\end{align*}
for all $(i,j)\in J$.

\vspace{5mm}

Out of the data given by the graph $G=(V,E)$ and the constant $k$ we compute the following integral constants
\begin{align*}
d_{G,k}&:= \bigg\lceil
\frac{(2k-3)\cdot |E|- |V|}{2}
\bigg\rceil
\intertext{and (with $n=|V|$) }
d_{n,k}'&:= \bigg\lceil
\frac{n(n-1)(2k-3)-n}{2}
\bigg\rceil
\end{align*}
(which we may incorporate in the codification of $(\Lambda V_{G,k},\dif)$).

Let us thus associate to $G$ and $k$, respectively to $(\Lambda V_{G,k},\dif)$ and $k$, yet another Sullivan algebra $(\Lambda W_{G,k},\dif)$, which extends $(\Lambda V_{G,k},\dif)$.
We set
\begin{align*}
W_{G,k}&:=V_{G,k}\oplus \langle z_1, \dots z_n\rangle
\intertext{with}
\deg z_i&=4(d_{n,k}'+n)+3
\intertext{and with}
\dif z_i&=x_i^{4(d_{n,k}'+n+1)} 
\end{align*}
for all $1\leq i \leq n$.

Obviously, this algebra can be constructed out of the algebra $(\Lambda V_{G,k},\dif)$ and the value $k$ in polynomial time.

\begin{rem}\label{rem01}
The finitely-generated algebra $(\Lambda W_{G,k},\dif)$ is pure and even elliptic, i.e.~its cohomology is finite-dimensional. This easily follows from the fact that, by construction, each form in $W_{G,k}^\even$---which necessarily defines a cohomology class, since $(\Lambda W_{G,k},\dif)$ is pure---represents a nilpotent cohomology class $[w]\in H(\Lambda W_{G,k},\dif)$. Indeed, we have $[x_i]^{4(d_{n,k}'+n+1)}=0$ for all $1\leq i\leq n$.
\end{rem}

We encode the graph $G$ as the number $n$ of its vertices together with the adjacency matrix representing the edges---an $(n\times n)$-matrix.

Recall that we encoded spaces by their minimal Sullivan models.

For fixed $k$, the spaces $|\langle (\Lambda V_{G,k},\dif) \rangle|$ and $|\langle (\Lambda W_{G,k},\dif) \rangle|$, i.e.~the spatial realisations (cf.~\cite{FHT01}.17) of the constructed minimal Sullivan algebras, thus have a codification the length of which is bounded by a polynomial in the length of the instance given by the graph. This means that our translations from graphs to algebras can be done in polynomial time. As we remarked, also the translation from $(\Lambda V_{G,k},\dif)$ to $(\Lambda W_{G,k},\dif)$ can be achieved with polynomial effort.

\begin{lemma}\label{lemma01}
The following assertions are equivalent
\begin{enumerate}[(i)]
\item
The Sullivan algebra $(\Lambda V_{G,k},\dif)$ is elliptic.
\item
The elements $[x_i]\in H(\Lambda V_{G,k},\dif)$ are nilpotent for all $1\leq i\leq n$.
\item
For all $1\leq i\leq n$ it holds that
\begin{align*}
[x_i]^{d_{G,k}+1}=0
\end{align*}
for all $1\leq i\leq n$.
\end{enumerate}
\end{lemma}
\begin{prf}

We shall prove that the ellipticity of the algebra is equivalent to the nilpotence of the $[x_i]$. The assertion on the order of the $[x_i]$ then can be deduced as follows:

If $(\Lambda V_{G,k},\dif)$ is elliptic, the formal dimension of $(\Lambda V_{G,k},\dif)$ is given by
\begin{align*}
&\sum_{(i,j)\in J} \deg y_{i,j}- \sum_{1\leq i\leq n} (\deg x_i-1)\\
=&  (2k-3)\cdot |E|- |V|\\
\leq & 2 d_{G,k}
\end{align*}
by \cite{FHT01}.32, p.~434. Consequently, by degree reasons, we obtain that $[x_i]^{d_{G,k}+1}=0$ for all $1\leq i\leq n$.

\vspace{5mm}

If $(\Lambda V_{G,k},\dif)$ is elliptic, the classes $[x_i]$ are necessarily nilpotent, since $H(\Lambda V_{G,k},\dif)$ is finite-dimensional by definition. Hence it only remains to prove the reverse implication in order to show the lemma.

Suppose that all the $x_i$ are nilpotent. We need to show that $H(\Lambda V_{G,k},\dif)$ is finite-dimensional. However, this algebra is finite-dimensional over $\qq$ if and only if it is finite-dimensional over $\cc$. Thus we may assume that we are using complex coefficients. Since the coefficient field then is algebraically closed and the algebra $(\Lambda V_{G,k},\dif)$ is simply-connected with $V$ finite-dimensional, we may use the criterion from proposition \cite{FHT01}.32.5, p.~439, saying that $(\Lambda V_{G,k},\dif)$ is elliptic if and only if every morphism
\begin{align*}
\varphi\co (\Lambda V_{G,k},\dif)\to (\cc[z],0)
\end{align*}
is trivial---here $\deg z=2$.

Such a morphism, however, is trivial on all degrees $V_{G,k}^{>2}=(V_{G,k}^{>2})^\odd$. Suppose it is not trivial in degree two, then it is given on a non-zero $x\in V^2_{G,k}$ by $x \mapsto \alpha z$ with $\alpha \in \cc$. Since the $[x_i]$ are nilpotent elements, there is a certain power of $[x]$ which vanishes. In other words, there exists an element $\tilde x$ in $\Lambda V$ with
\begin{align*}
0\neq \dif \tilde x\in \cc[x]\In \Lambda V_{G,k}
\end{align*}
This element has odd degree and $\varphi(\tilde x)=0$. Thus $\varphi$ does not commute with differentials; a contradiction.

\end{prf}

\begin{lemma}\label{lemma02}
If $[x_i]^{d_{G,k}+1}=0$ in $H(\Lambda W_{G,k},\dif)$ for all $1\leq i\leq n$, then we obtain an isomorphism of Sullivan algebras
\begin{align*}
(\Lambda W_{G,k},\dif) \cong (\Lambda V_{G,k},\dif) \otimes (\Lambda \langle z'_1, \dots , z'_n\rangle,0)
\end{align*}
with $\deg z'_i=4(d_{n,k}'+n)+3$.
\end{lemma}
\begin{prf}
Since $G$ is a simple (undirected) graph, there are at most $n \choose 2$ edges in $E$. We derive that
\begin{align*}
d_{G,k}&=\bigg\lceil\frac{(2k-3) |E|-n}{2} \bigg\rceil \\
&\leq \bigg\lceil\frac{(2k-3) \cdot n(n-1)/2-n}{2} \bigg\rceil\\
&=d_{n,k}'
\end{align*}
and infer the inequality
\begin{align*}
2(d_{G,k}+1) \leq 4(d'_{n,k}+n)+3
\end{align*}
(for $n\geq 1$).

Due to the fact that $(\Lambda W_{G,k},\dif)$ is identical to $(\Lambda V_{G,k},\dif)$ in degrees below degree $4(d'_{n,k}+n)+3$, we obtain that
\begin{align*}
H^{\leq 4(d'_{n,k}+n)+3}(\Lambda V_{G,k},\dif) \In H^{\leq 4(d'_{n,k}+n)+3}(\Lambda W_{G,k},\dif)
\end{align*}
is a graded subalgebra.

Consequently, since $[x_i]^{d_{G,k}+1}=0$ in $H(\Lambda W_{G,k},\dif)$ and since
\begin{align*}
\deg [x_i]^{d_{G,k}+1}=2(d_{G,k}+1)\leq 4(d'_{n,k}+n)+3
\end{align*}
we conclude that there is an element $\tilde x_i\in (\Lambda V_{G,k},\dif)$ with $\deg \tilde x_i=2 d_{G,k}+1$ and with $\dif \tilde x_i=x_i^{d_{G,k}+1}$ for each $1\leq i\leq n$.

It follows that for each $1\leq i\leq n$ the element
\begin{align*}
z_i':=z_i-\tilde x_i\cdot x_i^{2d'_{n,k}+2n-d_{G,k}+1} \in (\Lambda W_{G,k},\dif)
\end{align*}
is closed and not exact. The asserted splitting of differential graded algebras is a direct consequence.
\end{prf}

The main tool for proving theorem \ref{theoA} will be the following
\begin{prop}\label{prop01}
Suppose that $n,k\geq 1$. The following two statements are equivalent:
\begin{enumerate}[(i)]
\item
It holds that
\begin{align*}
\cupl (\Lambda W_{G,k},\dif)\leq d'_{n,k}+n
\end{align*}
\item
The Sullivan algebra $(V_{G,k},\dif)$ is elliptic.
\end{enumerate}
\end{prop}
\begin{prf}
We use the characterisation for the ellipticity of $(\Lambda V_{G,k},\dif)$ provided in lemma \ref{lemma01}.

If $(\Lambda V_{G,k},\dif)$ is elliptic, i.e.~if $[x_i]^{d_{G,k}+1}=0$ for all $[x_i]\in H(\Lambda V_{G,k},\dif)$ (with $1\leq i\leq n$), then lemma \ref{lemma02} yields the isomorphism
\begin{align*}
(\Lambda W_{G,k},\dif) \cong (\Lambda V_{G,k},\dif) \otimes (\Lambda \langle z'_1, \dots , z'_n\rangle,0)
\end{align*}
from which we directly see that the cup-length of $(\Lambda W,\dif)$ satisfies
\begin{align*}
\cupl(\Lambda W_{G,k},\dif)=\cupl(\Lambda V_{G,k},\dif)+n
\end{align*}
since $\deg z_i'$ is odd for $1\leq i\leq n$.

Since $(\Lambda V_{G,k},\dif)$ is elliptic with an element of minimal degree sitting in degree two, its cup-length can be estimated from above by its formal dimension divided by two, i.e.~in particular by
\begin{align*}
\cupl(\Lambda V_{G,k},\dif)\leq  2 d_{n,k}'/2=d_{n,k}'
\end{align*}
It follows that $\cupl(\Lambda W_{G,k},\dif))\leq d_{n,k}'+n$.

\vspace{5mm}

Conversely, if
\begin{align*}
\cupl(\Lambda W_{G,k},\dif)\leq d_{n,k}'+n
\end{align*}
then, in particular, the elements $[x_i]\in H(\Lambda W_{G,k},\dif)$ satisfy that
\begin{align*}
H(\Lambda W_{G,k},\dif) \ni [x_i]^{d'_{n,k}+n+1}=0
\end{align*}

Again we use that $(\Lambda W_{G,k},\dif)$ is identical to $(\Lambda V_{G,k},\dif)$ in degrees below degree $4(d'_{n,k}+n)+3$ and that
\begin{align*}
H^{\leq 4(d'_{n,k}+n)+3}(\Lambda V_{G,k},\dif) \In H^{\leq 4(d'_{n,k}+n)+3}(\Lambda W_{G,k},\dif)
\end{align*}
is a graded subalgebra, therefore.

Since
\begin{align*}
\deg x_i^{d_{n,k}'+n+1}=2(d_{n,k}'+n)+2\leq 4 (d'_{n,k}+n)+3
\end{align*}
for $n,k\geq 1$,
we conclude that
\begin{align*}
H(\Lambda V_{G,k},\dif) \ni [x_i]^{d_{n,k}'+n+1}=0
\end{align*}
(with $x_i$ now considered an element in $V_{G,k}$ and with the given power of its cohomology class already vanishing in $H(\Lambda V_{G,k},\dif)$).

Thus all the elements $[x_i]$ (for $1\leq i\leq n$) are nilpotent elements in $H(\Lambda V_{G,k},\dif)$. Due to lemma \ref{lemma01} it follows that the algebra $(\Lambda V_{G,k},\dif)$ is elliptic.
\end{prf}

The problem of $k$-colouring a graph, i.e.~attributing one of $k$ different colours to a vertex such that adjacent vertices have different colours,

\vspace{3mm}

\begin{itemize}
\item[$\mathcal{P}_2''$:]
Suppose that $k\geq 3$. Is the graph $G$ $k$-colourable?
\end{itemize}

\vspace{3mm}

is known to be $\mathbf{NP}$-complete (for $k\geq 3$)---for example cf.~\cite{GJ79}. In the proof of corollary \cite{LM00}.4, p.~92, it is shown that there is a polynomial reduction of $\mathcal{P}_2'$ to the problem

\vspace{3mm}

\begin{itemize}
\item[$\mathcal{P}_2'$:]
Given a simply-connected Sullivan algebra $(\Lambda V,\dif)$ with $\dim V<\infty$. Is it elliptic?
\end{itemize}

\vspace{3mm}

More precisely, the problem is reduced to

\vspace{3mm}

\begin{itemize}
\item[$\mathcal{P}_2$:]
Let $(\Lambda V_{G,k},\dif)$ be as constructed above. Does it constitute an elliptic algebra?
\end{itemize}

\vspace{3mm}

We are now ready to give the
\begin{proof}[\textsc{Proof of theorem \ref{theoA}}]
We consider the following decision problem.

\vspace{3mm}

\begin{itemize}
\item[$\mathcal{P}_3$:]
Let $(\Lambda W_{G,k},\dif)$ be an algebra as constructed above. Is the cup-length of $(\Lambda W_{G,k},\dif)$ smaller than or equal to $d'_{n,k}+n$?
\end{itemize}

\vspace{3mm}

There is a polynomial reduction of problem $\mathcal{P}_2''$ to problem $\mathcal{P}_2$. Problem $\mathcal{P}_2''$ is $\mathbf{NP}$-complete (for $k\geq 3$).
By proposition \ref{prop01} we can reduce problem $\mathcal{P}_2$ to problem $\mathcal{P}_3$ in polynomial time.
Hence we see that $\mathcal{P}_3$ is $\mathbf{NP}$-hard. However, the original problem $\mathcal{P''}$ is obviously harder than $\mathcal{P}_3$; thus it is $\mathbf{NP}$-hard.
\end{proof}

%%%%%%%%%%%%%%%%%%%%%%%%%%%%%%%%%%%%%%%%%%%%%%%%%%%%%%%%%%%%%%%%%%%%%%%%%%%%%%%%%%%

\section{Proof of theorem B}\label{sec03}

The proof of theorem \ref{theoB} will proceed along the lines of the proof of theorem \ref{theoA}.

\begin{prop}\label{prop02}
Suppose that $n,k\geq 2$. The following two statements are equivalent:
\begin{enumerate}[(i)]
\item
It holds that
\begin{align*}
\cat_0 (\Lambda W_{G,k},\dif)\leq d'_{n,k}+n
\end{align*}
\item
The Sullivan algebra $(V_{G,k},\dif)$ is elliptic.
\end{enumerate}
\end{prop}
\begin{prf}
Again, we use the characterisation for the ellipticity of $(\Lambda V_{G,k},\dif)$ provided in lemma \ref{lemma01}.

If $(\Lambda V_{G,k},\dif)$ is elliptic, i.e.~if $[x_i]^{d_{G,k}+1}=0$ for all $[x_i]\in H(\Lambda V_{G,k},\dif)$ (with $1\leq i\leq n$), then lemma \ref{lemma02} yields the isomorphism
\begin{align*}
(\Lambda W_{G,k},\dif) \cong (\Lambda V_{G,k},\dif) \otimes (\Lambda \langle z'_1, \dots , z'_n\rangle,0)
\end{align*}
Using theorem \cite{FHT01}.30.2.(ii) we compute the category of the tensor product as the sum of the categories
\begin{align*}
\cat(\Lambda W_{G,k},\dif)=\cat(\Lambda V_{G,k},\dif)+\cat(\Lambda \langle z'_1, \dots , z'_n\rangle,0)
\end{align*}
Since  $(\Lambda \langle z'_1, \dots , z'_n\rangle,0)$ is a formal algebra of cup-length $n$, its category also equals $n$ by example \cite{FHT01}.29.4, p.~388.

(Since this example is formulated for spaces, we observe that the rational category of a simply-connected space with rational homology of finite type is the category of its minimal Sullivan model---cf.~proposition \cite{FHT01}.29.4, p.~386.)

Since $(\Lambda V_{G,k},\dif)$ is elliptic with an element of minimal degree sitting in degree two, its Lusternik--Schnirelmann category can be estimated from above by its formal dimension divided by two---cf.~corollary \cite{FHT01}.29.1, p.385.

In particular, we obtain
\begin{align*}
\cat(\Lambda V_{G,k},\dif)\leq  2 d_{n,k}'/2=d_{n,k}'
\end{align*}
It follows that $\cat(\Lambda W_{G,k},\dif))\leq d_{n,k}'+n$.

\vspace{5mm}

Conversely, we assume that $\cat(\Lambda W_{G,k},\dif))\leq d_{n,k}'+n$. Either a straightforward direct check or a quote of proposition \cite{FHT01}.30.8.(ii), p.~410, for the trivial fibration with $F=X$ and $Y=\{\ast\}$ yields that
\begin{align*}
d_{n,k}'+n\geq \cat(\Lambda W_{G,k},\dif)\geq \cupl_0(\Lambda W_{G,k},\dif)
\end{align*}
Thus proposition \ref{prop01} yields that $(\Lambda W_{G,k},\dif))$ is elliptic.
\end{prf}

Hence we may establish the
\begin{proof}[\textsc{Proof of theorem \ref{theoB}}]
We consider the decision problem

\vspace{3mm}

\begin{itemize}
\item[$\mathcal{Q}_3$:]
Let $(\Lambda W_{G,k},\dif)$ be an algebra as constructed above. Is the Lusternik--Schnirelmann category of $(\Lambda W_{G,k},\dif)$ smaller than or equal to $d'_{n,k}+n$?
\end{itemize}

\vspace{3mm}

Again one uses the polynomial reduction of $\mathcal{P}_2''$ to $\mathcal{P}_2$ and the fact that $\mathcal{P}_2''$ is $\mathbf{NP}$-complete. Due to proposition \ref{prop02} we reduce $\mathcal{P}_2$ to $\mathcal{Q}_3$. Thus $\mathcal{Q}_3$ is $\mathbf{NP}$-hard. Again, the original problem $\mathcal{Q}''$ is harder than $\mathcal{Q}_3$, i.e.~it is $\mathbf{NP}$-hard, in particular.
\end{proof}

%%%%%%%%%%%%%%%%%%%%%%%%%%%%%%%%%% Bibliography %%%%%%%%%%%%%%%%%%%%%%%%%%%%%%%%%%%

%\bibliography{lib}{}
%\bibliographystyle{abbrv}

%\pagebreak\

\vfill

\begin{center}
\noindent
\begin{minipage}{\linewidth}
\small \noindent \textsc
{Manuel Amann} \\
\textsc{Department of Mathematics}\\
\textsc{University of Toronto}\\
\textsc{Earth Sciences 2146}\\
\textsc{Toronto, Ontario}\\
\textsc{M5S 2E4} \\
\textsc{Canada}\\
[1ex]
\textsf{mamann@uni-muenster.de}\\
\textsf{http://individual.utoronto.ca/mamann/}
\end{minipage}
\end{center}

\end{document}